**A Novel Highway Traffic Capacity Analyzing Method under Road Region Atmospheric Environment Constrains Based on Computational Fluid Dynamics Model**

**Ruohan Li**
School of Traffic and Transportation
Lanzhou Jiaotong University, Lanzhou, China, 730070
Email: rli04@villanova.edu

**Hualan Wang, Professor, Corresponding Author**
School of Traffic and Transportation
Lanzhou Jiaotong University, Lanzhou, China, 730070
Email: wanghualan126@126.com

**Qiyang Zhang**
School of Traffic and Transportation
Lanzhou Jiaotong University, Lanzhou, China, 730070
Email: 346632650@qq.com

**Ting Nie**
School of Traffic and Transportation
Lanzhou Jiaotong University, Lanzhou, China, 730070
Email: 1601239660@qq.com

## ABSTRACT

Highways always have a huge impact on the environment, quantifying the level of pollution and calculating the traffic capacity under environmental constraints is an important part of practicing environmental protection. Available traffic capacity methods do not focus on the traffic emissions impact on road region environment. To fill this research gap, this paper proposes a method consisting of Computational Fluid Dynamics (CFD) model and Traffic Flow Field and the Gas Flow Field coupled model for the traffic capacity model. The COPERT V model is adopted to simulate the traffic emissions on highways with field measurement data and traffic conditions as initial conditions. Then, these are served as the inputs, set and further calculated using coupled model in the CFD numerical simulation phase. Finally, modeling the Traffic Capacity model and calculating the traffic capacity. The experiment results in G30 Yongshan highway section demonstrate that, with the proposed method, the CFD simulation model for traffic emissions in highway region can reflect the practical conditions (the overall error is within 10%), while the traffic capacity can be well performance with the traffic capacity model under the atmospheric pollution constraints. The annual average traffic capacity of G30 Yongshan is 739,800 vehicles, and $SO_2$ under different constraints has the greatest limitation on traffic capacity (90,9662 vehicles), which is the most dominant impact in the environmental constraints of the road region, rest are ranked in: $NO_2$ (246,5373), $PM_{2.5}$ (323,6022), CO (771,8717), $CO_2$ (803,7217), and $PM_{10}$ (2202,0013).

## INTRODUCTION

Traffic pollution is one of the crucial problems of environmental pollution, and highway traffic pollutants become an important source of pollution in atmospheric environment pollution especially huge impact on highway region (*1–4*). However, quantifying the level of traffic pollution and calculating the traffic capacity under environmental constraints have always been the challenging task, for the complex and variable realistic environment and vehicle conditions (*5*) (The emissions are easily disturbed by factors including wind speed, wind direction, vehicle speed, driving behaviour, road conditions, fuel and traffic volume).

In recent years, researchers have turned their attention to traffic pollution, which is increasingly recognized as a serious problem. A great deal of research has been done on everything from soil heavy metal pollution (*6–9*), water pollution (*10*, *11*), noise pollution (*12–14*), Landscape ecology (*15*) to the familiar atmospheric pollution (*16–20*). In the study of traffic emissions on air pollution, the main focus has been on the concentration of traffic emissions in cities utilizing pollutant concentration prediction models includes the Gaussian diffusion model (*21*), and linear regression prediction model (*22*, *23*). Based on these concentration predictions, rough calculations of traffic emissions and the traffic capacity of the transportation environment in cities. It's worth noting that there have begun to study the traffic impact on the road area, usually, the study focuses on the impact of heavy metals (*9*, *10*) or vegetation (*24*), and is not able to analyze the combination of traffic and environmental pollution. And for the research method of traffic capacity of environment, e.g., Gray System Approach (*25*), system dynamics approach (*26*), hierarchical analysis (*27*), etc., which are biased towards traditional statistical analysis methods, and lack practical simulation methods to study the traffic carrying capacity intuitively and accurately. can't quantitatively study the relationship between the traffic volume and the distance of traffic pollution, the degree of environmental constraints.

For traffic capacity, they always focus on the capacity calculation of road resources, for more or better capacity allocation methods (*28–31*). But many people are committed to exploring the relationship between traffic emissions and environmental capacity. However, people seem to care only about the environmental goodness of the city they live in. Although some progress has been made in the cities, few have focused on the impact of highway traffic around the road region, which is often overlooked, especially for the ecologically functional or vulnerable areas.

Computational Fluid Dynamics (CFD) techniques used in Numerical simulation has shown in recent studies (*32–37*), creating a model of the surrounding environment and simulating traffic emissions under various conditions, makes it possible to understand the dispersion and distribution of traffic pollutants in different built environments. It has been widely used to explore environmental impacts of street or buildings in the cities (*5*, *38*, *39*). The basic idea, modeling the practical street or buildings model, the traffic pollution and other environmental factors are input as the boundary conditions, to simulate the traffic emission in a simulated scenario that resembles a real scenario. The CFD numerical simulation enables to provide a precise mathematical calculation, and better grasp the values of these pollutants at various points, which allows combining the traffic flow and pollution concentration for 1uantifying traffic capacity under environmental constraints.

In this paper, a CFD-based method was investigated for precise the traffic capacity calculation method under atmospheric environment constraints. The COPERT V model is adopted to simulate the traffic emissions on highways with field measurement data and traffic conditions as the initial conditions. Then, in the CFD numerical simulation phase where traffic emissions simulation on the highway region is performed, the traffic emissions and boundary conditions serve as the inputs, and the pollutant concentration constraints are set and further calculated using T-G coupled model calculating the traffic capacity under highway region atmosphere constraints. With this proposed method, the relationship between traffic volume and pollution concentration in the road region can be well captured.

The main contributions of this paper are as follows:

1. The Traffic Flow Field and the Gas Flow Field (T-G) coupled model are developed and adopted to further analyse capability of traffic volume under the atmospheric environment constraints. Without many extra parameters tuning, the T-G coupled model can easily change the boundary conditions in the CFD simulation model to adopt different constraints or environments;

2. This study proposes a traffic capacity method based on CFD numerical simulation and Traffic Flow Field and the Gas Flow Field (T-G) coupled model, and verified its effectiveness by comparison experiments, the simulation results and field measurement results can be in good agreement, with the error is within 10%;

3. The whole method is tested and verified taking G30 Yongshan section as the highway case, carrying out the annual traffic capacity under atmospheric environment constraints in G30 Yongshan. And analyzing from different constraints and vehicle types.

## METHODOLOGY AND MODELING

This study proposes a highway traffic capacity analyzing method using CFD numerical simulation under road region environment, which offers the maximum capability of highway traffic capacity under atmospheric environment constraints in the road region. The framework of the proposed method is illustrated in **Figure 1.** In the remaining parts of this section, each phase will be introduced in detail.

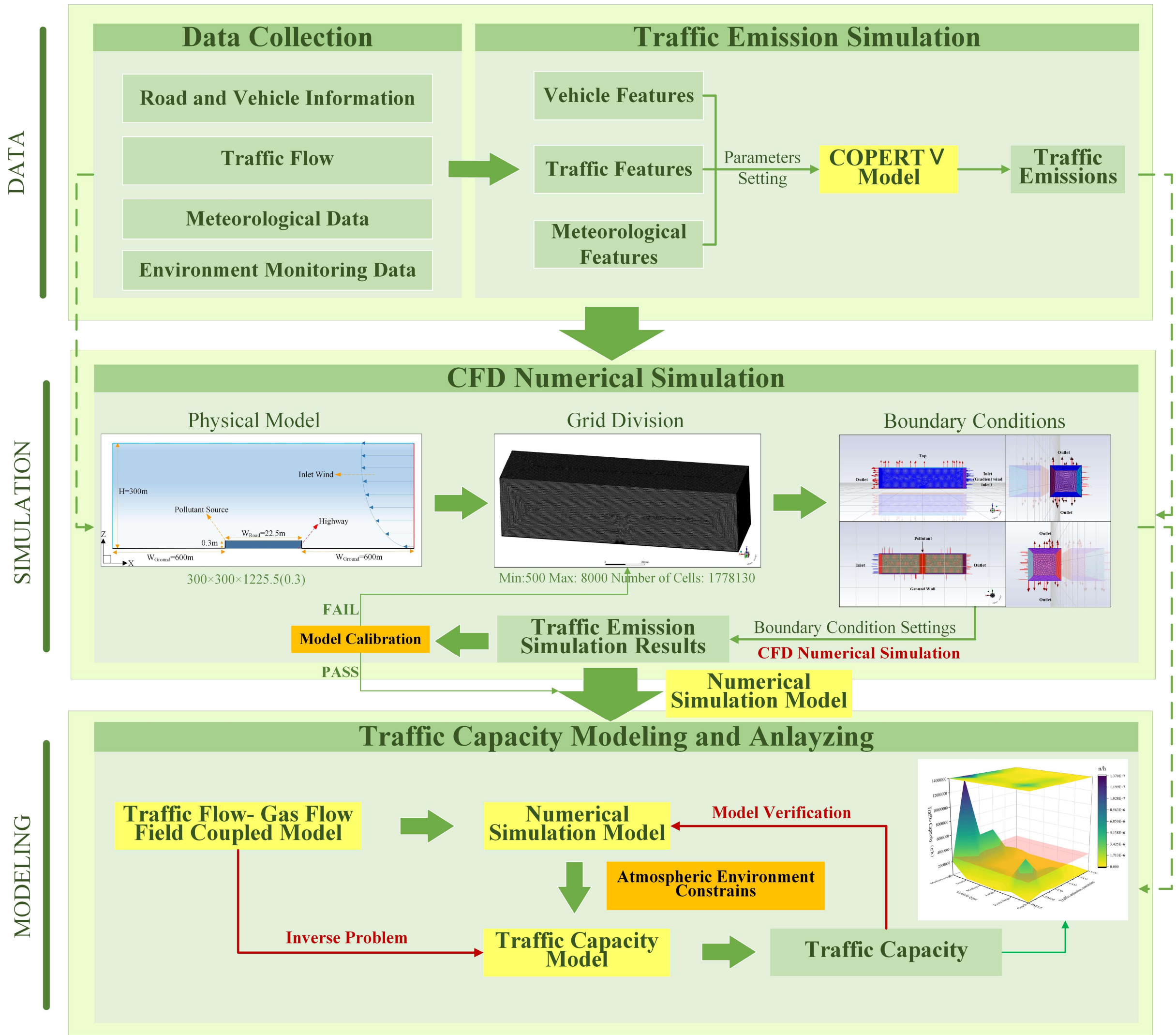

**Figure 1 The framework of the proposed method.**

### *A. Study Area and Data Collection*

**Study Area**: To better analyze the impact of highway traffic emissions on road region and to ensure that most of the pollutant data collected from traffic emissions, in this study, the G30-Yongshan (G30-YS) section of Lianhuo highway (located in Shandan County, Zhangye City, China) was selected. Shandan County is located in China's National Key Ecological Function Area (NKEFA) (**Figure 2**), where resource development should be restricted to maintain ecological functions, and therefore, traffic pollution should also be considered. Shandan is surrounded by mountains on three sides, and the climate is continental alpine and semi-arid, with an average temperature of 2°C~16°C. G30 is one of the main roads of China's highway network, which is located in Shandan County, and the G30-YS section with a total length of about 118 km, and its design speed of 100 km/h.

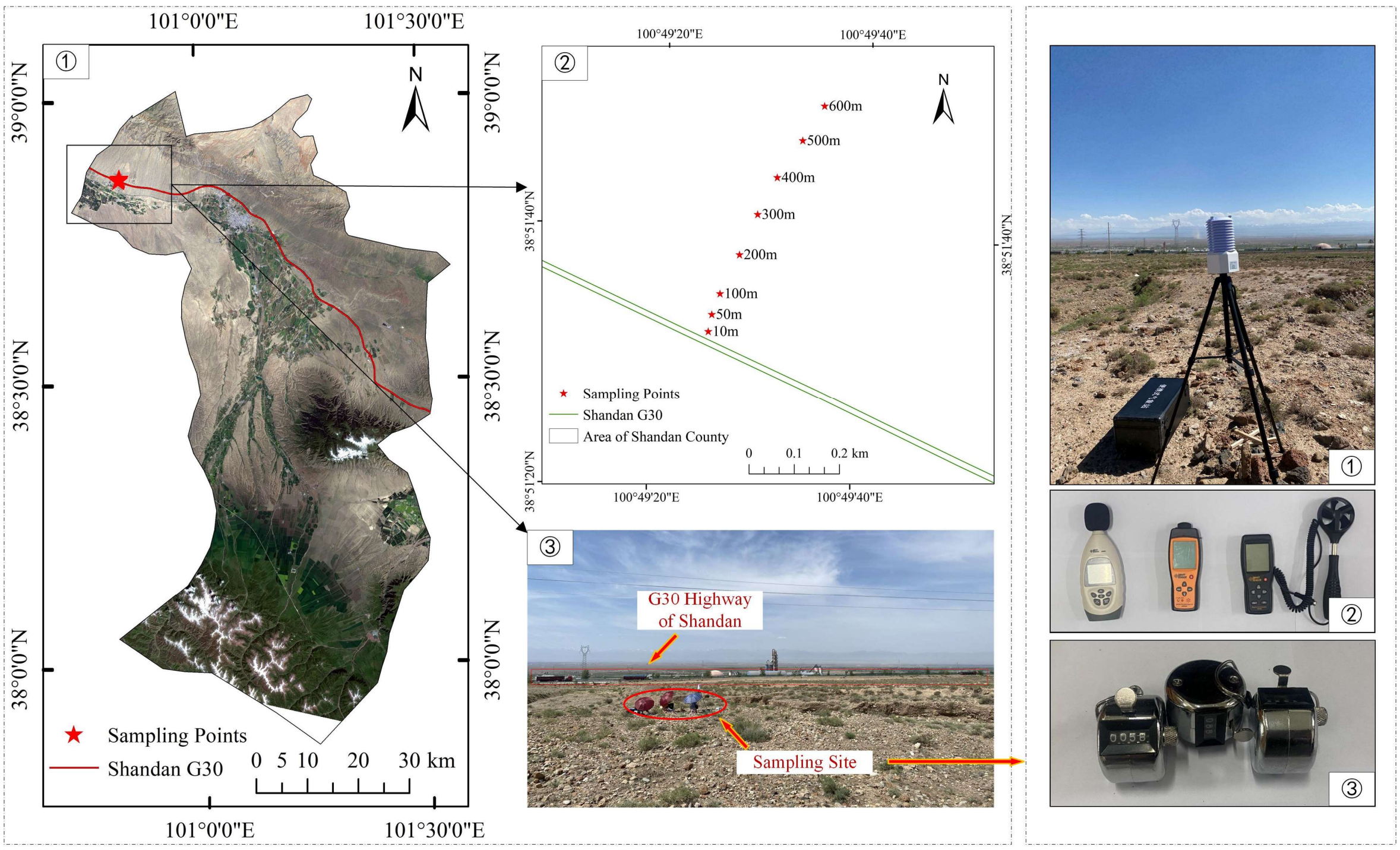

(a) (b)

**Figure 2 (a) The study area and sampling sites. (b) Field measurement equipment.** (a) ① The overview of G30 in Shandan county; ② The sampling sites in G30 road region; ③ The filed measurement in G30 road region. (b) ① Portable weather monitor; ② Noise detector, Carbon dioxide detector, Wind speed, and direction detector; ③ Traffic counter.

**Data Collection**: The data includes traffic volume data, and meteorological data (National Weather Science Data Center data). Other data were obtained from field measurements at the road region of G30-YS. The six conventional pollutants indexes: $PM_{2.5}$, $PM_{10}$, Carbon monoxide (CO), Ozone ($O_3$), Sulfur Dioxide ($SO_2$), and Nitrogen Dioxide ($NO_2$), where the pollutants indexes are continuously measured and recorded per minute pollutant values using a Tianhe portable weather monitor (TH-BQX6D).

For the field measurement, the portable weather monitor was used to collect traffic atmospheric pollutants, respectively, at 6 sampling sites (30m, 100m, 200m, 300m, 400m, 600m) from 30 m (100°49′28″*E*, 38°51′32″*N*) to 600m (100°49′37″*E*, 38°51′49″*N*) on one side of the G30-YS road region cross-section, the location of sampling sites in **Figure 2 (a)**. The pollutant monitors used here were calibrated before leaving the factory. Meteorological data, including temperature (°C), relative humidity (%), wind speed (m/s), wind direction (degree), air pressure (Pa), and six conventional pollutants were recorded at 1-min intervals by a TH-BQX6D portable weather monitor. Field measurements were taken over 3 days

(6/15/2022–6/17/2022) from 8:30 AM to 12:00 AM and 17:00 PM to 20:00 PM each day. This was a convenient time frame to observe the change in traffic flow and the concentration of traffic emissions. Before starting the official test, the research team conducted a preliminary survey on the daytime traffic volume of the highway and found that there was not much difference in the traffic volume (within the range of 2000 to 2500 vehicles per three-hour period) except for nighttime driving. Therefore, morning and afternoon were chosen as the measurement times.

### *B. Model Assumption*

To balance the computational cost with the simulation accuracy and simplify research, this study was based on several assumptions:

**Assumption 1:** Traffic Flow Fields carry out continuous traffic pollutant discharges in the form of surface sources he highway, and there is no loss of pollutant transfer between fields;

**Assumption 2:** Only highway traffic emissions, no other sources of pollution;

**Assumption 3:** Not taking into account factors such as the working conditions of the vehicle itself and the driving conditions of the driver;

**Assumption 4:** Air and emitted gaseous traffic pollutants are incompressible gas mixtures and obey the laws of diffusion and transport of gaseous flow fields in the atmosphere of the roadway; chemical changes occurring in the diffusion of traffic pollutants in the gaseous flow field are not taken into account; particulate matter is not deformed.

### *C. Traffic Emission Simulation*

**Traffic Conditions**: Traffic conditions for G30-YS from 2019 to 2021 were provided by Gansu Provincial Highway Administration. It is particularly noted that, in early 2020, the traffic volume from February to April was significantly lower than the traffic volume in 2019 and 2021 due to the impact of the COVID-19. In order to exclude this impact, the traffic volume in 2020 is taken as the average of the traffic volume in 2019 and 2021. The processed traffic volume data is shown in **Figure 3(a)**. In addition, G30-YS has a relatively fixed vehicle type, and the ratios of various vehicle types remain essentially the same between 2019 and 2021, so the vehicle type ratios are averaged over the three years, as shown in **Figure 3 (b)**. The average speed of G30-YS is taken as 90 km/h.

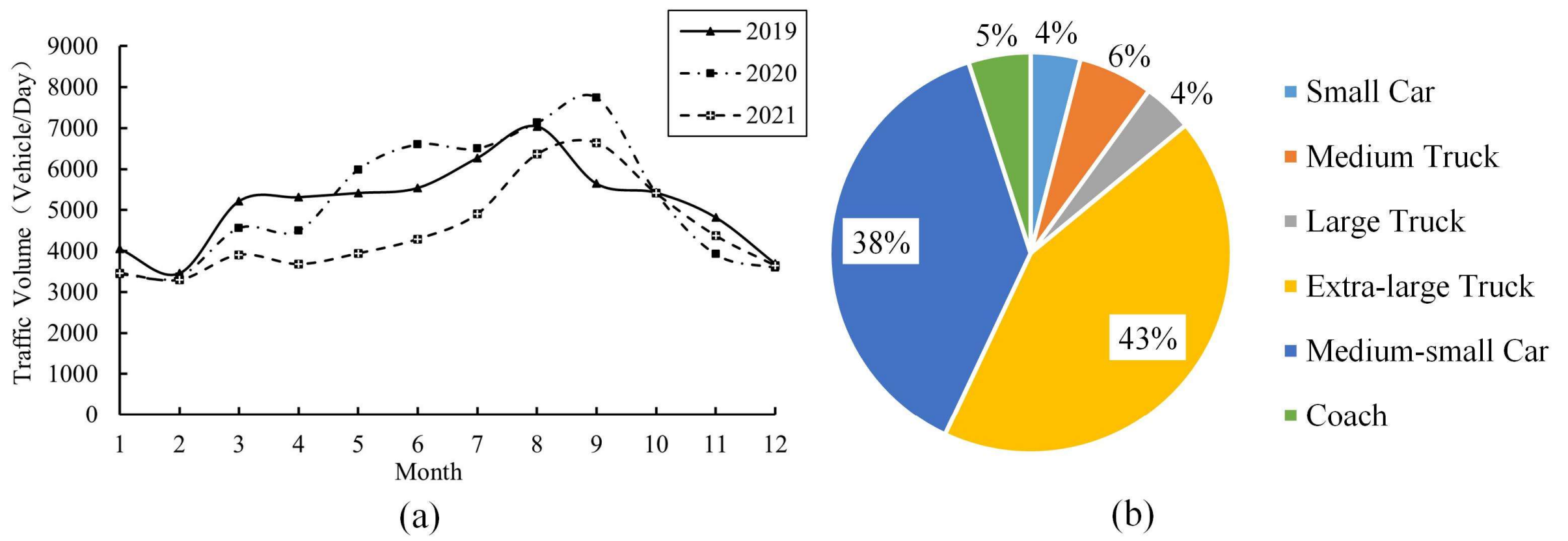

**Figure 3 (a) The traffic volume in G30-YS from 2019 to 2021. (b) Average vehicle type ratio in G30.**

**Simulation Setting**: Simulation of annual traffic emissions for G30-YS was carried out using COPERT 5.6. The file is created in COPERT with Asia-China as the background. Traffic emissions are calculated using the Tier 2 method, as more driving details are not considered in this study while simplifying the traffic capacity model. Minimum and maximum temperatures and humidity for each month of 2021 are added to the environmental information. The fuel in the COPERT model considers China's national standard for motor gasoline, National V (*Pollutant Emission Limits and Measurement Methods for Light-Duty Vehicles (China Stage V) GB18352.5-2013*). In the vehicle configuration, the vehicle types (Medium-small Car, Small Truck, Medium Truck, Large Truck, Extra-large Truck, and Coach) of the G30-YS are converted to the corresponding vehicle types in the COPERT model according to the correspondence, as shown in **TABEL 1**.

The average annual mileage of China's motor vehicles in medium-small car and bus have a certain trend every year, while the average annual mileage of other types of motor vehicles is more stable (*40*). Therefore, it can be used when the mileage travelled cannot be known exactly. In a survey of Chinese cities, Huo et al (*40*) estimated the average annual mileage to be 22,000 km in Beijing and 20,000 km in Tianjin, and Huang et al (*41*), in a study of motor vehicles in Shanghai, found that the average mileage of motor vehicles in Shanghai was 28,000km, and Huo et al (*42*) showed that it was 24,000 km in other cities. Therefore, with reference to the scholars' studies, the average annual mileage of road motor vehicles is set with the average annual mileage in **Table 1**.

**TABEL 1 Settings in COPERT**

| **Vehicle Type** | **Fuel** | **Category in COPERT** | **Euro Standard** | **Stock (n)** | **Average annual mileage (km)** |
|---|---|---|---|---|---|
| Medium-small Car | Petrol | Small | Euro 5 | 328361 | 14332 |
| Small Truck | Diesel | Rigid <=7,5 t | Euro V | 3986 | 30000 |
| Medium Truck | Diesel | Rigid 7,5 - 12 t | Euro V | 4971 | 35000 |
| Large Truck | Diesel | Rigid 20 - 26 t | Euro V | 1865 | 58000 |
| Extra-large Truck | Diesel | Rigid >32 t | Euro V | 119503 | 75000 |
| Coach | Diesel | Coaches Standard <=18 t | Euro V | 1752 | 60000 |

### *D. Numerical Simulation*

To investigate the impact of highway traffic emissions diffusion on the road region and to finally calculate the traffic capacity of the highway, the physical model is based on the field situation, taking the field sampling site as the center and the highway as the extension direction, and taking the highway to the left and right sides of the road for 150m each, which is a total length of 300m of the highway region simulation. Because of no buildings or plants in the vicinity of the measurement sites, and it's far from residential areas or factories, the highway region is relatively open and flat. In this paper, modeling traffic emission dispersion by treating the study area as a highway in a flat area.

**Turbulence model**: During the diffusion and transfer of traffic pollutants in the road region's atmospheric environment, the airflow always follows the physical conservation laws in nature, including conservation of mass, momentum, and energy. In addition, in FLUENT, it is a key step to choose the turbulence model reasonably. Turbulence models can be divided into three types: (1) one-equation turbulence model; (2) two-equation turbulence model; (3) Reynolds stress equation model. At present, the two-equation turbulence model plays an important role in the simulation calculation of fluids and practical engineering conditions. Turbulence modeling is mainly provided in FLUENT by the two-equation model *k*-epsilon (*k*-ε) model and k-omega (*k*-Ω) model. The *k*-ε model including Standard, RNG and Realizable, which converges more easily than the *k*-Ω model. Considering the specific experimental environment, the advantages of RNG *k*-ε cannot be fully utilized. Finally, the Standard *k*-ε turbulence model (**Equation 1**) was selected for the highway region model, which has been widely adopted in previous studies (*33*, *36*, *38*), to reveal the spatial distribution of wind field, turbulent flow, and particle concentrations in road region.

$$\frac{\partial}{\partial t}(\rho\varepsilon)+\frac{\partial}{\partial x_i}(\rho\varepsilon u_i)=\frac{\partial}{\partial x_j}\left[\left(\mu+\frac{\mu_t}{\sigma_\varepsilon}\right)\frac{\partial\varepsilon}{\partial x_j}\right]+C_{1\varepsilon}\frac{\varepsilon}{k}G_k-C_{2\varepsilon}\rho\frac{\varepsilon^2}{k} \quad (1)$$

where $\rho$ strands for the gas density, kg/m$^3$; $k$ is the turbulent kinetic energy, m$^2$/s$^2$; $G_k$ denotes the production term for the turbulent kinetic energy $k$ due to the mean velocity gradient; $\varepsilon$ is the dissipation rate, m$^2$/s$^3$; $C_{1\varepsilon}$, $C_{2\varepsilon}$ are constant; $\delta_k$, $\delta_\varepsilon$ is the turbulence Prandtl number for the $k$ and $\varepsilon$ equation; $\mu_t$ is the turbulent viscosity coefficient, $Pa \cdot s$; At high Reynolds numbers the equation is $\mu_t = \rho C_\mu k^2 / \varepsilon$, $C_\mu$ is constant.

**Models Descriptions**: The main consideration in the modeling is the atmospheric pollutants generated in the exhaust pipe of the vehicle, so the height of the exhaust pipe of the vehicle is considered. Referring to the heights used by Yu et al. (*34*) and Hang et al. (*38*) to study the traffic pollutant diffusion problem, and because the research road section has the largest proportion of extra-large trucks, and the height of the tailpipe of the large trucks is about 0.3m, and the pollutant emission is much larger than that of other models of vehicles such as minibuses, the height of the road pollutant source is set to 0.3 m. After a variety of onsite measurements, a 3D simulation model of the street canyon was developed with the CFD code FLUENT 6.3, shown in **Figure 4**.

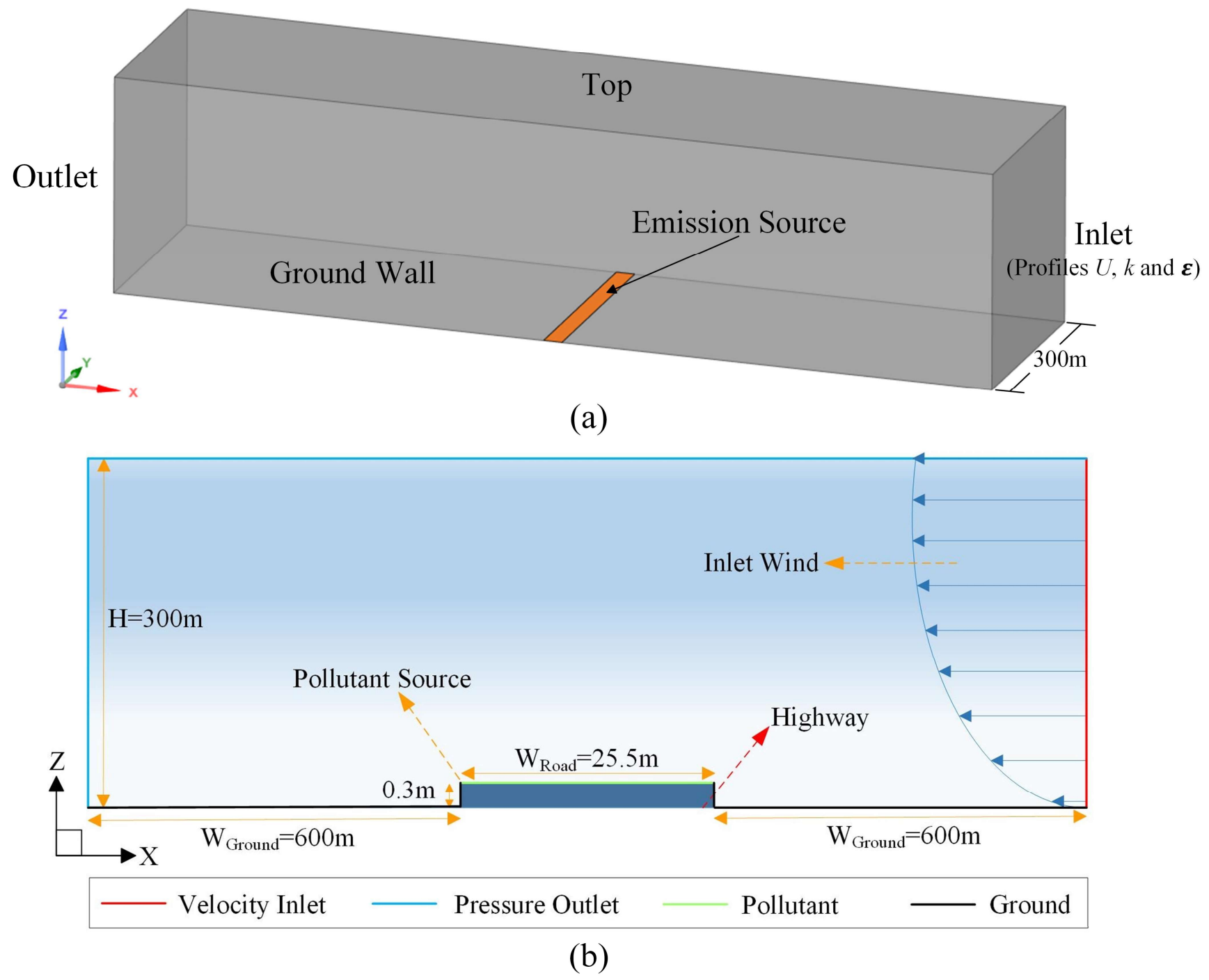

**Figure 4 Average vehicle type ratio in G30.**

**Grid Arrangement**: The mesh was created using a minimum mesh size of 500 mm for the tests and the maximum mesh size was selected from 500 to 15000 mm. To ensure grid convergence and numerical simulation accuracy, three final grid sizes were constructed: mesh1 (659234); mesh2 (1778130); and mesh3 (2989711). The mesh quality was evaluated using Skewness and Orthogonal quality, and computed under the same computer and working conditions to test for grid independence was tested.

**Figure 5** shows the results of the grid independence test for the highway centerline in the Y-axis direction (Y=0m, Z=2m, Y from 130m to 170m, for a total of 41 points). **Figure 5 (a)** and **(b)** show the velocity and CO concentration at the sampling points for the three grid sizes, respectively. The results are similar for the mesh2 and mesh3 grid sizes, where the maximum error in velocity and CO concentration is less than 5%, and this error is maximum for the merh1 grid size. A further increase in grid size from 1778130 to 659234 did not produce substantial changes in velocity and pollutant concentration, indicating the accuracy of the numerical solution and the independence of the network, and therefore a grid size of 1778130 (mesh2) was used.

Finally, taking into account the calculation speed and the accuracy of the calculation, the grid with the minimum neighbouring cell size of 0.5m, a maximum size of 5m, an expansion ratio of 1.2, and a total size of 1778130 is used for the meshing of the highway road area.

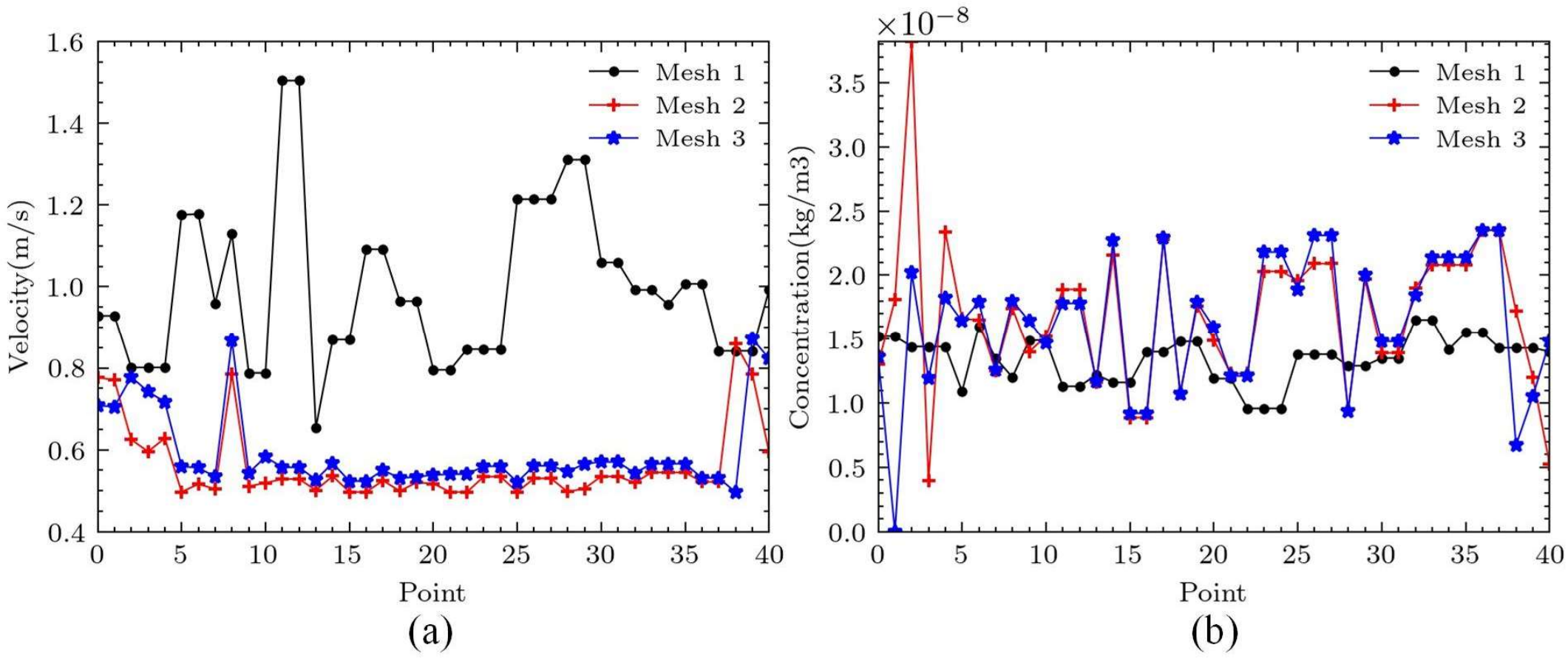

**Figure 5 Results of grid independence test. (a) velocity; (b) CO concentration.**

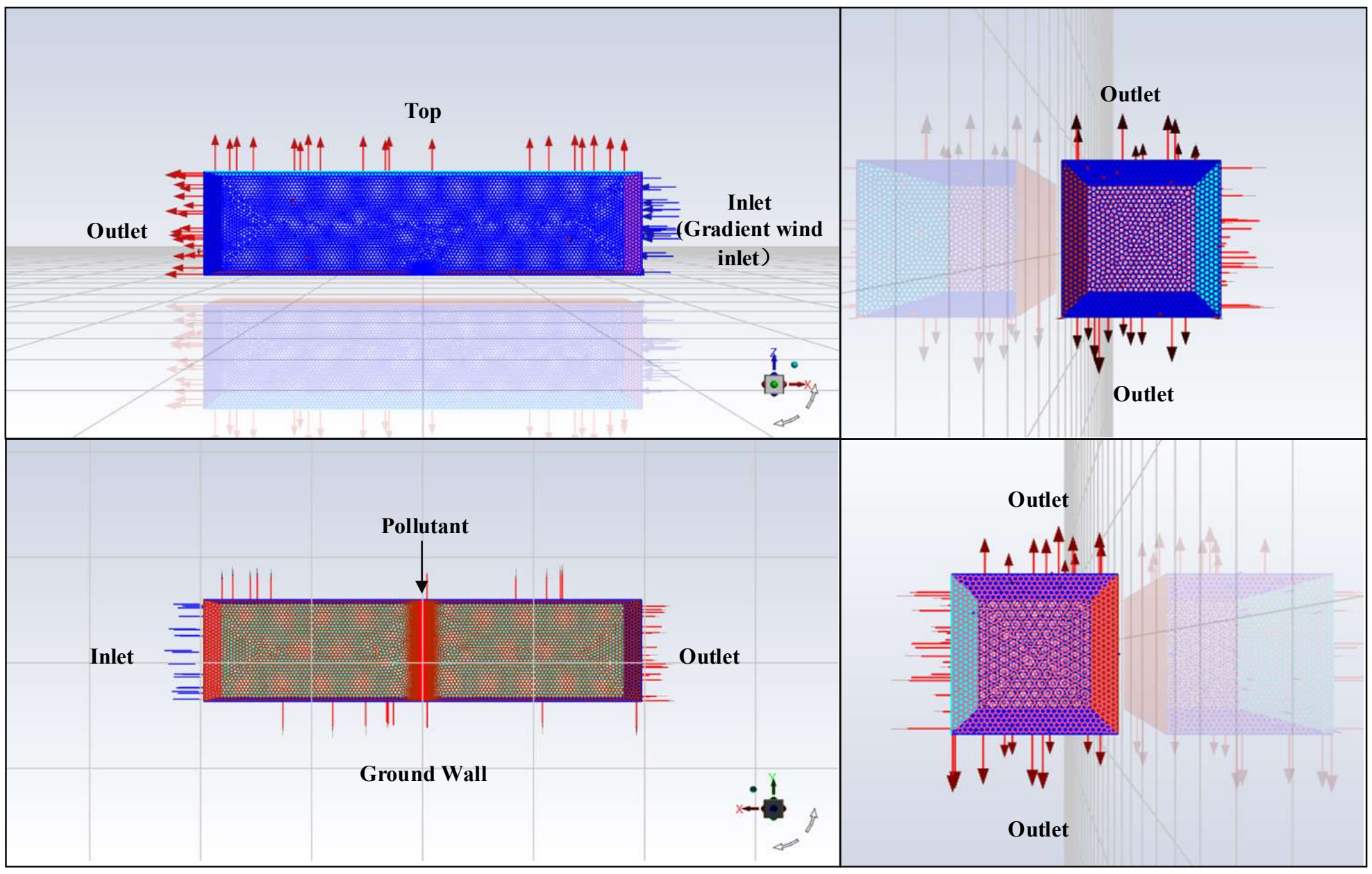

**Figure 6 Results of grid independence test. (a) velocity; (b) CO concentration.**

**Boundary Conditions**: The model boundary conditions mainly include the inlet and outlet settings and the traffic emission. The inlet of computational domain was set as velocity-inlet with medium turbulence intensity, the other three surfaces perpendicular to the vertical ground and top surfaces were set as pressure-outlet to approximate the actual atmospheric flow, and the rest of the surfaces in the region were set as wall surfaces (bottom and top) with a velocity boundary condition of no-slip (**Figure 6**). The pollutant emission inlet was "pollutant", and the particles were set as inert particles and subjected to gravity, CO, $CO_2$, $NO_2$, $PM_{2.5}$, and $PM_{10}$ were tested, and the diameters were set according to the diameter of the particles of each substance, which were set as $3.76\times10^{-10}$ m, $3.3\times10^{-10}$ m, respectively, $10.0\times10^{-10}$ m, $2.5\times10^{-6}$ m, and $1.0\times10^{-6}$ m. To restore highway traffic pollutant emissions as much as possible, the pollutant spraying direction was set to lay along the X-axis and Y-axis, and the spraying speed was set to 0.3 m/s by the whole pollutant surface source for uniform spraying.

**TABEL 1 Boundary conditions and CFD simulation settings**

**(a) Boundary conditions**

| **Boundaries** | **Concentration of pollutant** | **Velocity (m/s)** | ***k*(y)** | ***ε*(y)** |
|---|---|---|---|---|
| Inlet | 0 | $U(y)=\frac{u^*}{k}\ln(\frac{z+z_0}{z_0})$ | $K=\frac{U_*^2}{\sqrt{C_\mu}}$ | $\varepsilon=\frac{U_*^3}{\kappa(Z+Z_0)}$ |
| Top | $\frac{\partial C}{\partial n}=0$ | 0 | $\frac{\partial K}{\partial n}=0$ | $\frac{\partial \varepsilon}{\partial n}=0$ |
| Outlet | $\frac{\partial C}{\partial n}=0$ | 0 | $\frac{\partial K}{\partial n}=0$ | $\frac{\partial \varepsilon}{\partial n}=0$ |
| Ground Wall | $\frac{\partial C}{\partial n}=0$ | 0 | $\frac{\partial K}{\partial n}=0$ | $\frac{\partial \varepsilon}{\partial n}=0$ |
| Pollutant | 1 | 0.5 | | |

**(b) CFD simulation settings**

| **Settings** | **Specific content** |
|---|---|
| Computational domain | X×Y×Z; 1225.5m×300m×300m |
| Blockage ratio | <5% |
| Grid expansion ratio | 1.2 |
| Grid resolution | Minimum grid size: 0.5m<br>Maximum grid size: 5m<br>Grid number: 1778130 |
| Grid type | Unstructured grid |
| Gravity | -9.8N/kg |
| Inflow boundary conditions (Inlet) | Input wind velocity: 3.4m/s<br>Surface roughness factor α: 0.2<br>Turbulent kinetic energy $k$: 1.04$m^2s^{-2}$<br>Specific dissipation rate $\omega$: 1$s^{-1}$<br>Temperature: 32℃ |
| Outflow boundary (Outlet) | Zero gradient condition |
| Near wall treatment | STT k−$\omega$ model wall free function |
| Solving algorithms | SIMPLE algorithm for pressure-velocity coupling<br>Pressure using a physically weighted discrete scheme<br>Momentum k and $\omega$ adopt a first order upwind discretization |
| Relaxation factors | Default value in ANSYS Fluent |
| Convergence criteria | All variables: 1 $e^{-5}$ |

**(c) Settings of discrete phase injection source setup**

| Contents | CO | $CO_2$ | $NO_2$ | $PM_{2.5}$ | $PM_{10}$ | $SO_2$ |
|---|---|---|---|---|---|---|
| Pollutants(kg/s) | 8.02$e^{-04}$ | 0.7655 | 6.65$e^{-04}$ | 3.41$e^{-05}$ | 3.41$e^{-05}$ | 8.64$e^{-04}$ |
| Pollutant diameter | 3.76$e^{-10}$ | 3.3$e^{-10}$ | 10.0$e^{-10}$ | 2.5$e^{-6}$ | 1.0$e^{-6}$ | 2.46 $e^{-10}$ |

**Numerical Solver Settings**: The initial condition settings are set with reference to the field measurement meteorological conditions and local meteorological data including wind speed, wind direction, turbulence and pressure, details in **TABEL 1 (b)**.

The discrete phase is the pollutant emission factor, and six traffic emissions, $PM_{2.5}$, $PM_{10}$, CO, $CO_2$, $SO_2$ and $NO_2$, are considered for the discrete phase (DPM) in the road region environment. The DPM source is updated with an interval of 10 for each iteration, with a maximum number of steps of the tracking parameter of 50,000, and a step factor of 5. In the jet source, the pollutant surface is presented as a surface spraying. The diameter of the jet particles the total flow rate of release per second (kg/s), and the pollutant settings of the jet source are shown in the following **TABEL 1 (c)**. The material of the jet particles is selected according to the pollutant type, and its diameter distribution is uniform. The temperature of the jet source is set to 35°C, which is equal to the ground temperature. In the physical model, the trailing force criterion is the default spherical, and the particle agglomeration release method is also set with reference to the default value.

### *E. Traffic Capacity modeling*

The traffic capacity model is based on traffic emission simulation and CFD numerical simulation models, from the known traffic volume to calculate traffic emission concentration (positive problem) to deduce how to use the environmental pollution constraints inverse to seek the maximum capacity (inverse problem). Finally, established the traffic capacity model from the inverse problem, is used to solve a section of the highway environmental traffic capacity, and analyze the change of the capacity in different conditions (vehicle velocity, wind velocity, etc.).

**Coupled Model**: Viewing the diffusive transfer of highway traffic emissions to the road region atmosphere as a coupled problem. Then, the Traffic Flow Field and the Gas Flow Field are established separately. In Traffic Flow Field, vehicle operation on the highway is seen as a whole field, which reflects to the characteristics of the field itself and the sources of traffic pollutant emissions. Gas Flow Field includes the atmosphere on the highway and on the highway road region, which describes the diffusion and transfer of traffic pollutants in the road region atmosphere under the influence of external forces (wind, pressure, etc.) after the generation of traffic pollutants in the Traffic Flow Field. Therefore, The Traffic Flow Field can be regarded as the source field, and the process of pollutants generated by the Traffic Flow Field entering the Gas Flow Field can be regarded as a single coupling between the Traffic Flow Field and the Gas Flow Field. In this paper, the "Traffic Flow Field-Gas Flow Field" coupled model (T-G coupled model) is used to describe this single coupling to model the positive problem.

In the positive problem, the COPERT model was used as part of the Traffic Flow Field model. To reflect the characteristics of the field itself and the traffic pollutant emissions from the Traffic Flow Field, the establishment of the Traffic Flow Field is based on the traffic flow theory and basic theory of fields. The average traffic volume (Q), average vehicle speed (V) and average density (K) are taken as the basic variables of the Traffic Flow Field to give the unified expression of the traffic flow field in the highway road region, which is illustrated by **Equation 2**:

$$f(Q,V,K;W,I)=0 \quad (\Omega_T) \tag{2}$$

where $Q$ stands for the average highway traffic volume; $V$ is the average highway speed; $K$ is the Average highway density; $W$ is the set of output variable for Traffic Flow Field; $I$ is the set of input variables for the Traffic Flow Field; $\Omega_T$ is the definitional domain of the Traffic Flow Field equations.

Traffic emissions as output variables in the Traffic Flow Field:

$$W=\vartheta\cdot g_{\text{out}}(w_1,w_2,\ldots,w_n) \tag{3}$$

where $w_1,w_2,\ldots,w_n$ stands for the *n* traffic pollutants generated by the traffic flow field; $\vartheta$ is the traffic flow field coupling coefficient with other fields.

Traffic Flow Field can be illustrated by **Equation 4**:

$$f(Q,V,K;W,I)=\sum_{j}^{m} t_j\cdot Q_j\cdot\left(\mu_1\cdot\omega_{j,c}\cdot er_{i,j,c}+\mu_2 E_{CO_2,j}+\mu_3 E_{SO_2,j}\right) \tag{4}$$

where *t* stands for the travel time within a specific area; $\omega_{j,c}$ is the average annual miles traveled in category *j* on road *c*; $er_{i,j,c}$ is the emission factors for pollutants; $\mu_1=1$, when pollutant type *i* is CO, $NO_2$, $PM_{2.5}$ and $PM_{10}$, otherwise $\mu_1=0$; $\mu_2=1$, when pollutant type *i* is $CO_2$, otherwise $\mu_2=0$; $\mu_3=1$, when pollutant type $i$ is $SO_2$, otherwise $\mu_3 = 0$; $E_{CO_2,j}$ is $CO_2$ emissions from motor vehicles in category $j$; $E_{SO_2,j}$ is $CO_2$ emissions from motor vehicles in category $j$.

The controlling equations for the Gas Flow Field are those mentioned in **Governing Equations**, where $h(\rho,\gamma;o_{out},o_{in})$ is used to denote the Gas Flow Field.

Then, using Traffic Flow Field derived equations to describe pollutants generated by the Traffic Flow Field, and traffic pollutants transferred to the Gas Flow Field, the derivative equations of the Traffic Flow Field were illustrated by **Equation 5**:

$$q = \sum_{i}^{n} \mu_i w_i + \alpha w_{CO_2} + \beta w_{SO_2} \tag{5}$$

where $\mu_1$ is the control factor for traffic pollutant $i$; $w_i$, $w_{CO_2}$, $w_{SO_2}$ are emissions of transport pollutants in category $i$, emission of $CO_2$ and emission of $SO_2$; $\alpha$, $\beta$ are $CO_2$ emission factors and $SO_2$ emission factors, respectively.

The action of the Traffic Flow Field on the Gas Flow Field is a process of transfer of traffic pollutants, which is determined by the combination of the basic variables in the Traffic Flow Field and the output of the Traffic Flow Field, whose basic expression is shown in **Equation 6**:

$$C(W,o_{in}) = A \cdot f(Q,V,K;W,I) \tag{6}$$

where $C(W,o_{in})$ is the effect of the Traffic Flow Field on the Gas Flow Field; $A$ is the control coefficient of the action of the Traffic Flow Field on the Gas Flow Field, $A = (\mu_1,\dots,\mu_n,\alpha,\beta)$.

The final T-G coupled model can be illustrated by group of **Equation 1, 4, and 6**.

**Traffic Capacity Model**: The inverse problem is then modeled based on the positive problem model. The traffic volume (known condition) as in the positive problem is used as the variable to be solved in the inverse problem. The atmospheric environment constraints (pollutant concentration constrains) of the highway road region are added to the conditions of the inverse problem, then, the pollutant concentrations are known to solve for the traffic volume (unknown conditions), and the inverse problem model, the traffic capacity model is established.

The mathematical model of the Traffic Flow Field is expressed in the form of an inverse problem. The differential equations of the Gas Flow Field model are solved by the CFD numerical simulation, and the environment settings, boundary conditions, and input, output condition settings during the simulation process are used as the key parameter sets $E = (e_1, e_2, \dots, e_n)$ describing the Gas Flow Field in the traffic capacity model. The turbulence model boundary conditions and some other set numerical simulation parameters such as gravity, temperature, etc. in the CFD model are represented using the key parameter set as **Equation 7**:

$$E = (G, Temp, S_c, W_s, B) \tag{7}$$

where $G$ stands for is the gravitational force in the Gas Flow Field, N/kg; $Temp$ is the temperature in the gas flow field, ℃; $S_c$ is the turbulence model; $W_s$ stands for the wind speed of the gas flow field, m/s; $B$ is the boundary conditions.

In this paper, the constraints of road region atmospheric environmental pollutants $C_r$, are defined in accordance with *Atmospheric Environmental Quality Standard* (GB3095-2012) issued by China. The specific limit values are in accordance with the criteria for Class I, which includes I areas: nature reserves, scenic spots and other areas in need of special protection. Therefore, when maximum traffic capacity under the atmospheric environment constraints, traffic pollution concentrations should be less than $C_r$, using **Equation 8** to represent this relationship. When $\left[E \otimes c(W,o_m)\right] = C_r$, the traffic volume $T_{\max}$ in the Traffic Flow Field is the maximum traffic capacity.

$$\left[E \otimes c(W,o_m)\right] \le C_r \tag{8}$$

where $\left[E\otimes c\left(W,o_m\right)\right]$ stands for the CFD numerical simulation results of the coupled state of Gas Flow Field and Traffic Flow Field.

The traffic capacity finally can be illustrated by **Equation 9**:

$$\overline{T}_c=\sigma\cdot C_r\bigg/\left[E\otimes\sum_{j}^{m}t_j\left(90\cdot\mu_1\cdot er_j+\mu_2E_{j,co_2}+\mu_3E_{j,So_2}\right)\right] \tag{94}$$

where $\overline{T}_c$ stands for the average traffic capacity; $\sigma$ is the conversion factor, $\sigma=31.536$.

## EXPERIMENTS AND RESULTS

### *A. Model Calibration*

To verify the accuracy of the numerical simulation results, numerical simulation results were verified to ensure that the model established by numerical simulation is as accurate as possible. The numerical simulation concentration sampling was carried out with the field measurement point arrangement. The numerical simulation results were compared with the field measurement, shown in **TABEL 2**.

**TABEL 2 Comparison of simulation and field measurement results**

| Traffic emissions | Distance | Background | Field measurement(μg/m$^3$) | Concentration(μg/m$^3$) |
|---|---|---|---|---|
| CO | 30m | 0.16 (mg/m$^3$) | 0.348 | 0.101 |
| | 100m | | 0.437 | 0.093 |
| | 200m | | 0.312 | 0.199 |
| | 300m | | 0.676 | 0.169 |
| | 400m | | 0.404 | 0.192 |
| | 600m | | 0.352 | 0.094 |
| $CO_2$ | 30m | 638.23 (mg/kg) | 410.647 | 175.849 |
| | 100m | | 423.452 | 192.747 |
| | 200m | | 435.588 | 234.004 |
| | 300m | | 444.882 | 239.360 |
| | 400m | | 450.000 | 227.325 |
| | 600m | | 443.176 | 225.687 |
| $NO_2$ | 30m | 6.68 (μg/m$^3$) | 0.000893 | 0.151 |
| | 100m | | 0.001010 | 0.155 |
| | 200m | | 0.002339 | 0.199 |
| | 300m | | 0.001754 | 0.193 |
| | 400m | | 0.001169 | 0.191 |
| | 600m | | 0.014035 | 0.143 |
| $SO_2$ | 30m | 2.37 (μg/m$^3$) | 2139.057 | 915.234 |
| | 100m | | 1983.965 | 334.081 |
| | 200m | | 690.863 | 156.524 |
| | 300m | | 445.133 | 21.669 |
| | 400m | | 40.283 | 20.042 |
| | 600m | | 398.807 | 16.521 |
| $PM_{2.5}$ | 30m | 52.66 (μg/m$^3$) | 19.482 | 6.923 |
| | 100m | | 19.393 | 10.175 |
| | 200m | | 15.626 | 5.812 |
| | 300m | | 15.725 | 7.498 |
| | 400m | | 15.918 | 6.060 |
| | 600m | | 15.491 | 8.108 |
| $PM_{10}$ | 30m | 78.56 (μg/m$^3$) | 24.982 | 4.698 |
| | 100m | | 24.444 | 3.514 |
| | 200m | | 19.702 | 7.956 |
| | 300m | | 19.842 | 6.370 |
| | 400m | | 20.041 | 7.226 |
| | 600m | | 19.544 | 3.295 |

In **TABEL 2,** $CO_2$, $NO_2$, $PM_{2.5}$ and $PM_{10}$ in the background data are all greater than in the actual detection, this is because that the background concentration in the environment is not added in the simulation. However, CO and $SO_2$ are higher than the meteorological monitoring surface values in the actual detection, the meteorological monitoring of CO is 0.16 mg/m$^3$, $SO_2$ is 2.37 mg/m$^3$, while the field measurement of CO average value is 0.42 mg/m$^3$, $SO_2$ average value is 949.68 mg/m$^3$, exactly they are typical of the traffic emissions, which proves that the traffic pollution emissions as the main source of pollution in the study area.

Finally, based on the numerical simulation of **TABEL 2**, the background concentration (70% of the atmospheric primary standard) was added. **Figure 7** shows that the numerical simulation results after adding the background values remain basically the same as the actual collected data, especially the results of $CO_2$ and $PM_{2.5}$ are basically in agreement with the measured results, with an error within 1%, and the values of the remaining three pollutants are also within the acceptable range (the overall error is within 10%), which shows the effectiveness of the CFD numerical model.

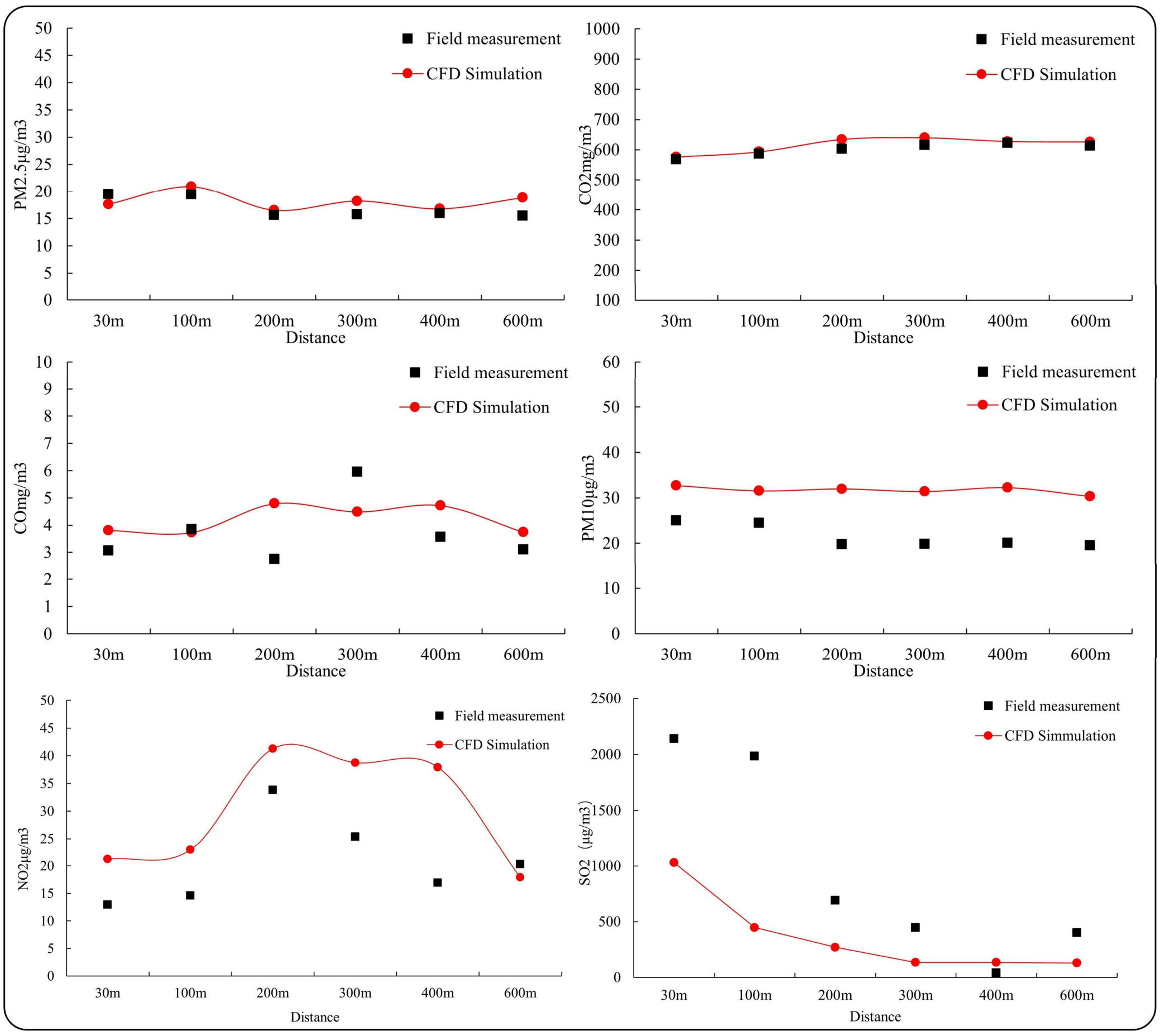

**Figure 7 Comparison of field measurement and CFD simulation results**

### *B. Calculation Results*

To calculate the traffic capacity for one year, the 24-hour concentration limit value was used as the pollutant concentration constraint, and a distance of 600 m from the highway was used as the environmentally constrained comparison point. Using **Equation 9** to calculate the traffic volume when the CFD simulation value at the 600m point is equal to the environmental constraint value to get the annual traffic capacity, illustrated in **TABEL 3**.

**TABEL 3 Results.**

**(a) Pollutant concentration values at different distances**

| Site | Distance | CO (μg/m3) | $CO_2$ (μg/m3) | $NO_2$ (mg/m3) | $PM_{2.5}$ (mg/m3) | $PM_{10}$ (μg/m3) | $SO_2$ (μg/m3) |
|---|---|---|---|---|---|---|---|
| CO Constraint | | | | | | | |
| point-30 | 30 | 0.4297 | 748.29 | 0.6425 | 16293 | 11056. | 1151.1 |
| point-100 | 100 | 0.3957 | 820.20 | 0.6595 | 23946 | 8270.2 | 657.54 |
| point-200 | 200 | 0.8468 | 995.76 | 0.8468 | 13678 | 18724 | 313.16 |
| point-300 | 300 | 0.7191 | 1018.55 | 0.8212 | 17646 | 14991 | 158.38 |
| point-400 | 400 | 0.8170 | 967.34 | 0.8127 | 14262 | 17006 | 151.16 |
| point-600 | 600 | 0.4297 | 748.29 | 0.6425 | 16293 | 11056 | 1151.11 |
| Constraint | 600 | 4 | 1500 | 40 | 35 | 40 | 20 |
| $CO_2$ Constraint | | | | | | | |
| point-30 | 30 | 0.4475 | 779.17 | 0.6691 | 30.6752 | 20.8164 | 1198.6148 |
| point-100 | 100 | 0.4120 | 854.04 | 0.6868 | 45.0845 | 15.5702 | 684.6790 |
| point-200 | 200 | 0.8817 | 1036.85 | 0.8818 | 25.7524 | 35.2523 | 326.0850 |
| point-300 | 300 | 0.7488 | 1060.58 | 0.8552 | 33.2230 | 28.2249 | 164.9220 |
| point-400 | 400 | 0.8507 | 1007.25 | 0.8463 | 26.8513 | 32.0177 | 157.3993 |
| point-600 | 600 | 0.4165 | 999.99 | 0.6336 | 35.9258 | 14.5998 | 141.1239 |
| Constraint | 600 | 4 | 1500 | 40 | 35 | 40 | 20 |
| $NO_2$ Constraint | | | | | | | |
| point-30 | 30 | 0.1373 | 239.01 | 0.2052 | 9.4093 | 6.3852 | 200.22 |
| point-100 | 100 | 0.1264 | 261.97 | 0.2107 | 13.8292 | 4.7760 | 114.37 |
| point-200 | 200 | 0.2705 | 318.05 | 0.2705 | 7.8993 | 10.8133 | 54.47 |
| point-300 | 300 | 0.2297 | 325.33 | 0.2623 | 10.1908 | 8.6577 | 27.55 |
| point-400 | 400 | 0.2610 | 308.97 | 0.2596 | 8.2363 | 9.8211 | 26.29 |
| point-600 | 600 | 0.1278 | 306.74 | 0.1944 | 11.0198 | 4.4783 | 23.57 |
| Constraint | 600 | 4 | 1500 | 40 | 35 | 40 | 20 |
| $SO_2$ Constraint | | | | | | | |
| point-30 | 30 | 0.0507 | 88.1875 | 0.0757 | 3.4717 | 2.3559 | 135.6610 |
| point-100 | 100 | 0.0466 | 96.6617 | 0.0777 | 5.1025 | 1.7622 | 77.4930 |
| point-200 | 200 | 0.0998 | 117.3519 | 0.0998 | 2.9145 | 3.9897 | 36.9068 |
| point-300 | 300 | 0.0848 | 120.0379 | 0.0968 | 3.7600 | 3.1944 | 18.6661 |
| point-400 | 400 | 0.0963 | 114.0024 | 0.0958 | 3.0389 | 3.6236 | 17.8147 |
| point-600 | 600 | 0.0471 | 113.1810 | 0.0717 | 4.0659 | 1.6523 | 15.9726 |
| Constraint | 600 | 4 | 1500 | 40 | 35 | 40 | 20 |
| $PM_{2.5}$ Constraint | | | | | | | |
| point-30 | 30 | 0.1326 | 230.817 | 0.1982 | 9.0873 | 6.1667 | 624.838 |
| point-100 | 100 | 0.1221 | 252.998 | 0.2035 | 13.3559 | 4.6125 | 356.923 |
| point-200 | 200 | 0.2612 | 307.151 | 0.2612 | 7.6289 | 10.4432 | 169.988 |
| point-300 | 300 | 0.2218 | 314.181 | 0.2533 | 9.8420 | 8.3614 | 85.974 |
| point-400 | 400 | 0.2520 | 298.384 | 0.2507 | 7.9545 | 9.4850 | 82.052 |
| point-600 | 600 | 0.1234 | 296.234 | 0.1877 | 10.6427 | 4.3251 | 73.568 |
| Constraint | 600 | 4 | 1500 | 40 | 35 | 40 | 20 |
| $PM_{10}$ Constraint | | | | | | | |
| point-30 | 30 | 1.2261 | 2134.7 | 1.8331 | 84.0426 | 57.0319 | 3283.91 |
| point-100 | 100 | 1.1290 | 2339.9 | 1.8816 | 123.5206 | 42.6586 | 1875.85 |
| point-200 | 200 | 2.4158 | 2840.7 | 2.4158 | 70.5555 | 96.5828 | 893.39 |
| point-300 | 300 | 2.0516 | 2905.7 | 2.3429 | 91.0229 | 77.3294 | 451.85 |
| point-400 | 400 | 2.3308 | 2759.6 | 2.3187 | 73.5661 | 87.7209 | 431.24 |
| point-600 | 600 | 1.1411 | 2739.8 | 1.7360 | 98.4280 | 40.0000 | 386.65 |
| Constraint | 600 | 4 | 1500 | 40 | 35 | 40 | 20 |

**(b) Average annual traffic capacity**

| Vehicle type | Type ratio | $PM_{2.5}$ Constraint | $PM_{10}$ Constraint | CO Constraint | $CO_2$ Constraint | $NO_2$ Constraint | $SO_2$ Constraint |
|---|---|---|---|---|---|---|---|
| Medium-small car | 62.00% | 2006334 | 13652408 | 4785605 | 4983075 | 1528532 | 563990 |
| Small truck | 6.50% | 210341 | 1431301 | 501717 | 522419 | 160249 | 59128 |
| Medium truck | 6.00% | 194161 | 1321201 | 463123 | 482233 | 147922 | 54580 |
| Large truck | 3.50% | 113261 | 770700 | 270155 | 281303 | 86288 | 31838 |
| Extra-large truck | 19.50% | 631024 | 4293903 | 1505150 | 1567258 | 480748 | 177384 |
| Coach | 2.50% | 80901 | 550500 | 192968 | 200930 | 61634 | 22742 |
| The sum of a year (n) | | 3236022 | 22020013 | 7718717 | 8037217 | 2465373 | 909662 |
| Equivalents of standard vehicle (n) | | 4748862 | 32314370 | 11327219 | 11794619 | 3617935 | 1334929 |

### *C. Discussion*

Different types of vehicles are converted to standard vehicles to get the number of equivalent standard vehicles. **TABEL 3 (b)** shows the annual average traffic capacity under different pollutant constraints, the maximum annual average traffic capacity is 32314370 vehicles under $PM_{10}$ constraint, and the minimum annual average traffic capacity is 1334929 vehicles under $SO_2$ constraint. The current annual average traffic volume of the G30-YS section of the highway is 2661896 vehicles, and only the annual average traffic capacity under the $SO_2$ constraint is lower than the actual annual average traffic volume, while under the remaining pollutant constraint, the annual average traffic capacity of the rest of the pollutant constraint, except for $NO_2$, is much larger than the actual value. This indicates that the limitations of the annual average traffic capacity of the road area are mainly in the two pollutants $SO_2$ and $NO_2$.

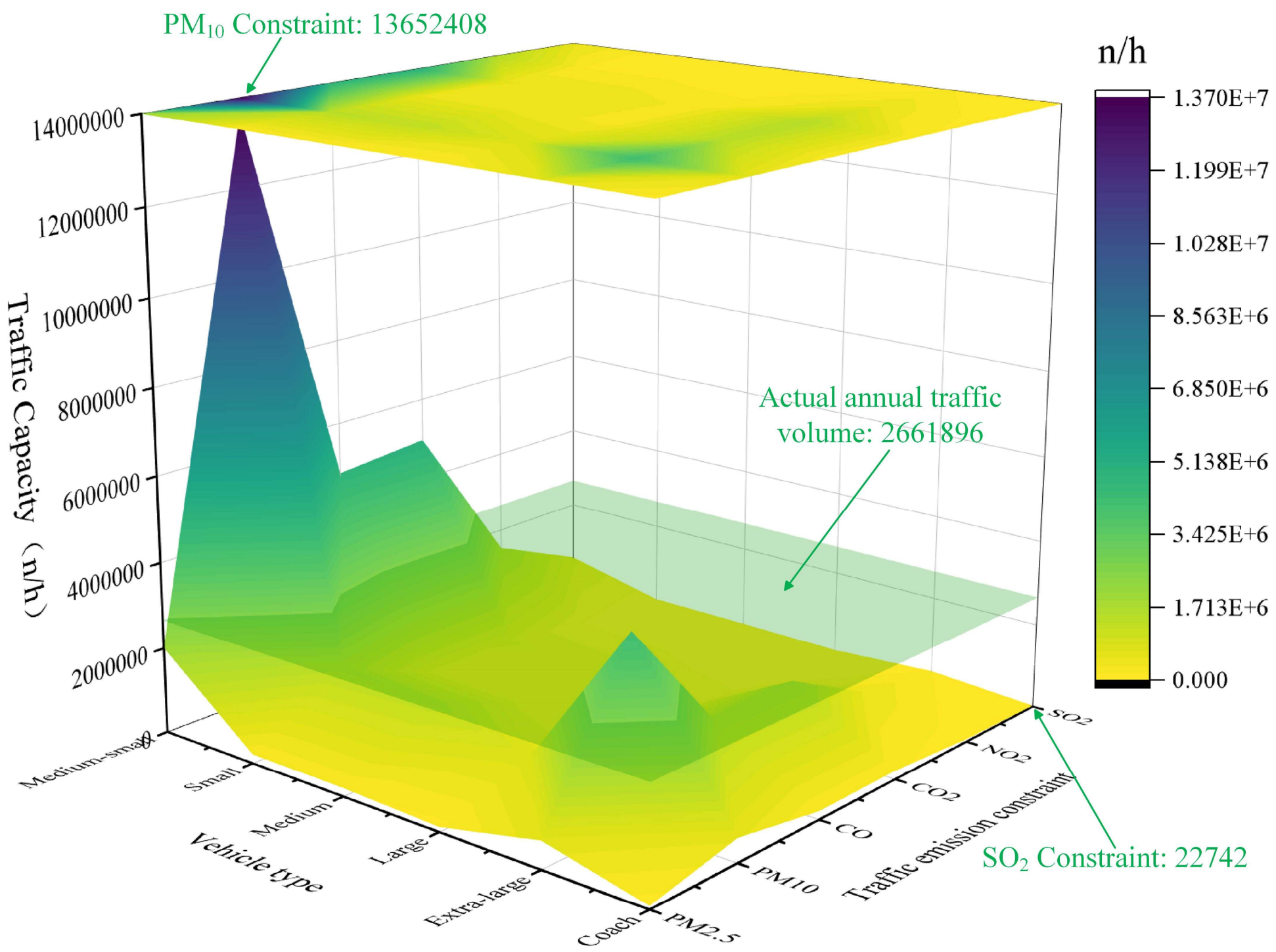

**Figure 8 Average annual traffic capacity in G30-YS**

Among all the pollutant constraints, $SO_2$ has the greatest limitation on the traffic capacity. In China, it is common for motor vehicles to operate under inadequate combustion conditions, mainly due to the low speed of motor vehicles, poor operating conditions, engines often working in a fuel-rich state, high CO pollutant emissions, and high CO concentrations. Such a situation is more mainly in urban areas, urban roads are narrow and congested, road running cars are often in low-speed operation or idling state. At the same time, there are many buildings in the city, and the roads become "canyons" in the midst of these buildings, making it difficult for pollutants to diffuse. The CO constraint of the study section is not the main limitation, probably because of the high speed of vehicles running on the highway, the engine fuel can be fully combusted, and the highway around the open and unobstructed, and the wind speed, highway traffic emissions of pollutants can get a better diffusion.

$SO_2$ became the most significant limitation, mainly because of the large percentage of trucks and the high percentage of diesel fuel burned in this highway section. Vehicle emissions of $SO_2$ are fuel-related. Generally speaking, diesel engines emit more sulfur dioxide than gasoline engines. Sulfur dioxide in the air in contact with water will form "acid rain", seriously endangering human health. For trucks that are mainly burning diesel-based, compared with gasoline engines, diesel engines have low fuel consumption, high energy efficiency, high reliability, low-speed performance, torque, tools, such as the use of a wide range of advantages. The high emission area of diesel vehicles has a high degree of overlap with the serious regional pollution zones in China. Therefore, there is no exception on G30-YS, where $SO_2$ becomes the main constraint factor, then followed by $PM_{2.5}$ and $PM_{10}$.

## CONCLUSIONS

In this paper, a novel highway traffic capacity analyzing method with CFD model is proposed for analyzing the maximum traffic volume under environment capacity. With the proposed method by modeling the CFD numerical simulation model and traffic capacity mode, it is very easy to simulate the maximum traffic capacity of a flat and open highway under environmental constraints and to analyze the change of traffic carrying capacity under different conditions by changing the boundary conditions in the CFD model. Using the G30-YS section as an example study, we performed field measurement and simulation modeling, calculated the traffic capacity of the area under the respective pollutant constrains. The results show that, the simulated values of the CFD model can well match the pollutant concentrations measured in the field, and it is feasible and effective to use the CFD modeling method to calculate the traffic capacity. The results of the calculations can show the traffic capacity in terms of different vehicle types and different pollutant constraints.

It is witnessed that when simulating some pollutant emissions in CFD model, i.e., SO2 and $NO_2$, the model might not match well with field measurements. In practice, different pollutant concentration values may be generated due to different wind speeds, wind directions or vehicle operating conditions. To tackle this problem and further enhance the model robustness, for future studies it is suggested to investigate domain generalization and increase sample and study cases to more comprehensive modeling.